\documentclass[11pt]{article}
\usepackage{multirow}
\usepackage{amsfonts, amssymb, graphicx, a4}
\usepackage{amsmath,dsfont}
 \usepackage[usenames, dvipsnames,svgnames]{xcolor}
 \usepackage{hyperref}
\def\inst#1{$^{#1}$}

\newtheorem{theorem}{Theorem}[section]
\newtheorem{lemma}[theorem]{Lemma}
\newtheorem{proposition}[theorem]{Proposition}

\newtheorem{definition}[theorem]{Definition}
\newtheorem{corollary}[theorem]{Corollary}
\newtheorem{remark}[theorem]{Remark}
\newtheorem{remarks}[theorem]{Remarks}

 \newenvironment{proof}{\noindent {\bf Proof.\,}
 }{\hspace*{\fill}$\qed$\medskip}

\def \P {{\mathbb P}}
\def \E {{\mathbb E}}

\newcommand{\cX}{{\cal X}}

\newcommand{\cvd}{\hfill\raisebox{3pt}
{\fbox{\rule{0mm}{1mm}\hspace*{1mm}\rule{0mm}{1mm}}\,}\vspace{8pt}}

\def \a {{\alpha}}

\def \l {{\lambda}}
\def \m {{\mu}}

\def \g {{\gamma}}
\def \t {{\tau}}

\def \d {{\delta}}
\def \p {{\pi}}

\def\eqref#1{(\ref{#1})}

\newcommand{\be}[1]{\begin{equation}\label{#1}}
\newcommand{\ee}{\end{equation}}

\newcommand{\bl}[1]{\begin{lemma}\label{#1}}
\newcommand{\el}{\end{lemma}}

\newcommand{\br}[1]{\begin{remark}\label{#1}}
\newcommand{\er}{\end{remark}}

\newcommand{\bt}[1]{\begin{theorem}\label{#1}}
\newcommand{\et}{\end{theorem}}

\newcommand{\bd}[1]{\begin{definition}\label{#1}}
\newcommand{\ed}{\end{definition}}

\newcommand{\bcl}[1]{\begin{claim}\label{#1}}
\newcommand{\ecl}{\end{claim}}

\newcommand{\bp}[1]{\begin{proposition}\label{#1}}
\newcommand{\ep}{\end{proposition}}

\newcommand{\bc}[1]{\begin{corollary}\label{#1}}
\newcommand{\ec}{\end{corollary}}

\newcommand{\bi}{\begin{itemize}}
\newcommand{\ei}{\end{itemize}}

\newcommand{\ben}{\begin{enumerate}}
\newcommand{\een}{\end{enumerate}}
 \newcommand{\bpr}{\begin{proof}}
 \newcommand{\epr}{\end{proof}}

\def \qed {{\square\hfill}}
\def \Z {{\mathbb Z}}

\def \P {{\mathbb P}}
\def \E {{\mathbb E}}

\def\1{\rlap{\mbox{\small\rm 1}}\kern.15em 1}

\def\build#1_#2^#3{\mathrel{\mathop{\kern 0pt#1}\limits_{#2}^{#3}}}
\def\tend#1#2#3{\build\hbox to 12mm{\rightarrowfill}_{#1\rightarrow
#2}^{#3}}

{}
\def\converge#1#2#3{\build\hbox to
15mm{\rightarrowfill}_{\hbox{\scriptsize #3}}^{#1\rightarrow #2}}
\def\converg#1#2#3{\build\hbox to
15mm{\rightarrowfill}_{\hbox{\scriptsize #3}}^{#1\uparrow #2}}

\begin{document}

\title{Exact results on the first hitting \\ via strong metastability times. }

%\begin{aug}
\author{F.\ Manzo\inst{1}\and
E.\ Scoppola\inst{2}}

\maketitle
\begin{center}
{\footnotesize
\vspace{0.3cm} \inst{1}  Supported by Dipartimento di Matematica e Fisica,
 University of Roma
``Roma Tre''\\
Largo San Murialdo, 1 - 00146 Roma, Italy\\

\vspace{0.3cm} \inst{2} Dipartimento di Matematica e Fisica, University of Roma
``Roma Tre''\\
Largo San Murialdo, 1 - 00146 Roma, Italy\\
\texttt{scoppola@mat.uniroma3.it}\\ }

\end{center}

\begin{abstract}
In the setting of non-reversible Markov chains on
finite or countable state space,
exact results on the
distribution of the  first hitting time to a given set $G$ are obtained. A new notion of  ``strong metastability time"
is introduced to describe the 
 local relaxation time. 
This time is defined via
a generalization of the  strong stationary time to a   ``conditionally strong quasi-stationary time"(CSQST).

Rarity of the target set $G$ is not required and the initial distribution can be completely general. 
The results  clarify the the role played by the initial distribution on the exponential law; they  
are used to give a general notion of metastability and to discuss
  the relation between the 
exponential distribution of the first hitting time and metastability.

\end{abstract}

%\tableofcontents
{\bf Keywords:} first hitting, strong stationary time, metastability time.

\section{Introduction }

The first hitting problem  is
 a key ingredient in the discussion of metastability in the framework of Markov processes.
The distribution of   the first hitting time $\t_G$ to a goal set $G$,    is widely discussed in the literature. 
In particular (see for instance \cite{A}, \cite{Kei79}, \cite{AB1},  \cite{AB2}) it is well known that under suitable ``rarity hypotheses" for the target set $G$ and for a suitable starting measure,
the hitting time $\tau_G$ is approximately exponential.  

In this paper we study the distribution of the hitting time in the general setting of non-reversible Markov chains, starting from an arbitrary initial distribution $\a$ and we prove an exact (non-asymptotic) representation formula for $\mathbb{P}\left({\tau^\a_G} > t\right)$ in terms of a new notion of  ``strong metastability time".
This representation formula provides an explicit, exact, probabilistic interpretation for the error terms in the exponential approximation and
it gives a new control on the role of the starting distribution.
We use the notion of strong metastability time to give a general definition of metastability in terms of the counterposition between a ``global relaxation time" and a ``local relaxation time".

Heuristically, in the metastable case, the system, before reaching $G$, thermalizes to a local  equilibrium.  From then on, $G$ is reached after many trials, which give the exponential behavior of $\tau_G$.
This means that,  for the exponential behavior of $\tau_G$, it is sufficient  that the time needed to reach $G$ is  much longer than the time needed to relax to the local equilibrium. 

This time comparison is very common in the literature, see for instance \cite{BodH}.  Rarity hypotheses and, more generally,  metastability hypotheses are often  given in terms of the ratio between two different time-scales:
a ``short" time-scale  characterizing  the approach to some local equilibrium, often described by the  {\it{quasi-stationary measure}}
\begin{equation}
\mu^{*}(\cdot):=\lim_{t\to\infty}\mathbb{P}\left(X_{t}=\cdot \ | \ t<\tau_{G}\right)\label{mu*},
\end{equation}
and a ``long" time-scale
characterizing the arrival to $G$.
The precise definition of the short and long time-scales, however, vary according to the methods used by different authors and to the different regimes at issue. 
 
In some of these regimes, hitting times are a very powerful tool to describe the behavior of the chain (see e.g. \cite{OV}).
This is the case when the invariant measure $\pi$ and $\mu^*$ asymptotically concentrate on single points: if $\pi $ concentrates in $G$ and $\mu^* $ concentrates in a point $m$, we can identify $G$ with the stable state and $m$ with the metastable state.  
In this case, we can state the metastability hypothesis (see \cite{FMNS} for a discussion) in terms of a time comparison as
\be{rar0}
  \sup_{x\neq m,G } \frac{\mathbb{E}_{x} \tau_{m \cup G}}{\mathbb{E}_{m} \tau_{ G}} \longrightarrow 0.
\ee
The idea is that if we observe the system on a time scale larger than the local relaxation
time $R:= \sup_{x\neq m,G }\mathbb{E}_{x} \tau_{m \cup G}$,
 the process behaves like a two state chain, since all other points $x\not \in\{m,G\} $ decay rapidly to $m\cup G$.

The notion of local relaxation time that we give in this paper is strictly inspired 
by the time $\tau_{m\cup G}$, but, unlike all other choices in the literature, it leads to an exact
representation formula for the hitting time in the general case:
non reversible, non recurrent, non asymptotic, for any initial state.

The key idea is to replace the hitting time to the metastable state $m$ with  a sort of ``hitting time to the 
{quasi-stationary measure}" obtained via  a generalization of  the notion of strong stationary time. 
We  define  a {\it conditionally-strong quasi-stationary time} (CSQST in the following) $ \t_*$
satisfying 
$$
 \P(X^{\alpha}_t=y,\, \t^{\alpha}_*=t)=\m^*(y)\P(\t^{\alpha}_*=t<\t^{\alpha}_G) \qquad \forall \;y\not\in G, \quad \forall t\ge 0.
$$
CSQST's are the central object of this paper and give a very powerful description of the approach to the local equilibrium.
  
The point is  to use this CSQST in the decomposition
\be{basici}
 \P(\t^{\alpha}_G>t)=
 \P(\t^{\alpha}_G>t \ ; \ \t^{\alpha}_{*}\le t)
 +\P( {\t^{\alpha}_{*,G}}> t)
\ee
 where 
  $\t_{*,G}:=\t^{\alpha}_G\wedge \t^{\alpha}_*$, called   {\it strong metastability time},  takes the role of
  $\t_{m\cup G}$.
  Equation (\ref{basici}),  
 has some interesting features that we will exploit in this paper in order to obtain bounds on the exponential approximation:
 \begin{enumerate}
 	\item It is an exact formula, that does not require reversibility and does not need any approximation or asymptotic.  
 	\item The conditional quasi-stationary property of $\tau_*$ allows to give exponential bounds on the first term
	in the r.h.s. of (\ref{basici}).  Since $\lambda^t = \P(\t^{\mu^*}_G>t )$, and the event in the first term implies a visit to a measure proportional to $\mu^*$,  its probability can be evaluated in terms of $\lambda^t$. 
	\item The second term $\P( \t^{\alpha}_{*,G}> t)$  gives a probabilistic interpretation of the error term in the exponential approximation. 
	Exponential behavior emerges when this second term is negligible with respect to the first one; therefore, it is natural to express the metastable hypothesis in terms of the comparison between the
	mean values 
	of  $\t^{\alpha}_{*,G}$ and of  $\t^{\alpha}_{G}$.  The first represents the local
	relaxation time,  to be compared with the hitting time.

 	\item The role of the starting measure $\alpha$ is  explicit.  We will show that, in the long period, the initial state gives a time-shift:  some states help and some hinder to reach $G$. To our knowledge, no other result in the literature gives a comparable control on the initial state.  
 	In many approaches based on renewal ideas, the lack of control  on the effect of the starting state is a primary source of error propagation in the exponential approximation.
 \end{enumerate}

Let us mention that the idea of using a strong time that somehow
catches the arrival to the quasi-stationary measure is not new in
the literature; in \cite{DM09}, for a birth-and-death process starting from $0$, in a particular
regime, the authors construct what they call a \emph{``strong
quasi-stationary time''} for this purpose. 
Although the motivations are similar, our approach is different, our notion of Conditionally Strong Quasi Stationary Time is completely
general and its existence does not require any additional assumptions
besides ergodicity of the stochastic matrix outside $G$.
\bigskip

The paper is organized as follows.
\bi
\item[--] In Section \ref{localch} we introduce a local chain $\tilde X_t$ on $A:= G^c$ related
to the Doob transform of $P$. This construction is useful in order to
\bi
\item[-]
determine the dependence of $\P(\t^\a_G>t)$ on the initial distribution $\a$ in terms of a time shift
$\d_\a$;
\item[-]
control the distribution of the CSQST $\t^\a_*$ and the terms in (\ref{basici}). In particular
 its first term  can be rewritten as $\P(\t^{\alpha}_G>t \ ; \ \t^{\alpha}_{*}\le t) = \lambda^{t+\delta_\alpha} \left( 1 - \tilde{s}^{\tilde \alpha} (t) \right)$ , where $\delta_\alpha$ is the time-shift that depending on  $\alpha$,  $ \tilde{s}^{\tilde \alpha} (t)$ is a separation from stationarity for the local Markov chain $\tilde X_t$  and $\tilde \alpha$ is the measure on $A$ induced by $\alpha$. 	
\item[-] obtain rought estimates on $\P(\t^{\alpha}_G>t )$ in terms of $\tilde{s}^{\tilde \alpha} (t)$.
\ei
\item[--] In Section \ref{proofth} we collect the proofs of our main results.
We introduce an auxiliary process, the {\it tracking process}, to provide a
 construction of the CSQST, which is discussed in Section \ref{s_main1}.

 In Section \ref{s_main2} we prove the representation formula for the hitting time.

 We use again the  tracking process to construct  the {\it ephemeral measure} in Section
\ref{ephemeral} describing
the process before the CSQST. Even if the tracking process is not Markovian, the ephemeral measure,
constructed with it, has a nice semigroup property that turns out to be the main ingredient in the
proof of submultiplicativity of the distribution of $\t^\a_{*,G}$, 
the local relaxation time.
\item[--] We give in Section \ref{ex} a simple example where the CSQST is explicitly constructed
in terms of a sequence of hitting times. This example is also useful to discuss the relation
between metastability and exponential distribution of the decay time.
\ei

\bigskip

  \bigskip
   
%%%%%%%%%%%%
%%%%%%%%%%%%%%%%%%
\subsection{General setting, definitions and preliminary remarks}
\label{1.1}
\subsubsection{Notation}
\label{s_not}
We collect in this subsection definitions and  notation used in the paper.
\begin{itemize}
\item
{\bf Process:}
we will consider  a discrete time Markov chain  $\{X_t\}_{t\in \mathbb N}$  on a finite state space $\cX$.
Our results can be extended to the case of countable state space but for the sake of simplicity we consider the finite case.
We  denote  by $P(x,y)$ the
transition matrix
  and  by $\m^x_t(\cdot)$ the  measure at time $t$, starting at $x$, i.e.,  
$\m^x_t(y)\equiv\P(X^x_t=y)\equiv P^t(x,y)$, for any $y\in \cX$.  More generally, given an initial distribution
${\a}$ on $\cX$
$$
\m^{\a}_t(y)=\P(X^{\a}_t=y)=\sum_{x\in\cX}{\a}(x)P^t(x,y)
$$
Starting conditions (starting state $x$ or starting measure $\a$) will be denoted by a  
superscript in random variables  (i.e., $X^x_t$, $X^\a_t$,  $\t^x$, ...).

Let $G\subset \mathcal{X}$ be a target set and  $\tau_{G}$   its first hitting time 
$$\tau_{G}:=\min\{t\ge0\:;\: X_{t}\in G\}.$$

\item
{\bf Separation:}
given two measures $\nu_1$ and $\nu_2$ on $\cX$  their
\emph{separation} is defined by
\be{dsep}
sep(\nu_1,\nu_2):= \max_{y\in\cX}\Big[1-\frac{\nu_1(y)}{\nu_2(y)} \Big]
\ee
\item
{\bf Scalar product:}
given two functions $a(x)$ and $b(x)$ on $A=\cX\backslash G$ we define their {\it scalar product} as
$$
a\cdot b:=\sum_{x\in A}a(x)b(x).
$$
\item
{\bf Strong Stationary Time} (see \cite{AD1} and \cite{AD2}):
a randomized stopping time $\tau^\alpha_\pi$ is a \emph{Strong Stationary Time (SST)}   for the
Markov chain $X^\a_t$ with starting distribution $\a$ and stationary measure $\p$,  if for any $t\ge 0$ and
$y\in\cX$
\[
\mathbb{P}\left(X^{\alpha}_t=y,\tau^{\alpha}_\pi=t\right)=\pi(y)\mathbb{P}\left(\tau^{\alpha}_\pi=t\right).
\]
This is equivalent to say
$$
\mathbb{P}\left(X^{\alpha}_t=y\big|\tau^{\alpha}_\pi\le t\right)=\pi(y)
$$
If $\tau^\alpha_\pi$ is a strong stationary time  then
\be{minsst}
 \P(\tau^{\alpha}_\pi>t)\ge sep(\m^\a_t,\pi),\qquad \forall t\ge 0
\ee
When the identity holds in (\ref{minsst})  the strong stationary time is {\it minimal}.
\item
{\bf Ergodicity:}
We will study   the process  $\{X_t\}_{t\in \mathbb N}$ up to time  $\tau_{G}$,  so it is not restrictive to consider $G$ as a set of absorbing states.
We assume ergodicity on $A:=\mathcal{X}\backslash G$. More precisely,
 denoting by $[P]_A$  the sub-stochastic matrix obtained by $P$ by restriction to $A$
$$
[P]_A(x,y)=P(x,y) \qquad \forall x,y\in A,\qquad\hbox{ with }\quad \sum_{y\in A}[P]_A(x,y)\le 1,
$$
we suppose $[P]_A$ a {\it primitive matrix}, i.e., there exists an integer $n$ such
that $\big([P]_A\big)^n$ has strictly positive entries.
\item
{\bf  Quasi-stationary measure on $A$:}
by the Perron-Frobenius theorem,   there exists $\l<1$ such that $\l$ is the spectral radius of
  $[P]_A$ and there exists a unique non negative left eigenvector of $[P]_A$ corresponding to $\l$, i.e.,
  \be{eqmqs}
  \m^*[P]_A=\l \m^*
  \ee
   this is called {\it quasi-stationary measure}.
 We get immediately 
 $$\P\left(\tau^{\mu^*}_{G}>t\right)=\l^t.$$
Moreover,
  the {quasi-stationary measure} $\m^*$ satisfies the following equation  (see \cite{DS}): 
\begin{equation}
\mu^{*}(\cdot)=\lim_{t\to\infty}\mathbb{P}\left(X^x_{t}=\cdot \ | \ t<\tau^x_{G}\right)\qquad\forall x\in A.\label{mu*'}
\end{equation}
\item
{\bf Hitting distribution:} starting from $\m^*$, the {\it hitting distribution} to $G$ is defined as
\be{defomega}
\omega(y):=\P\Big(X^{\m^*}_1=y\Big| X^{\m^*}_1\in G\Big)= \frac{\sum_{z\in A}\m^*(z)P(z,y)}{1-\l}.
\ee

\item
{\bf Conditionally-strong quasi-stationary time:}
a randomized stopping time $\tau^\alpha_*$ is a \emph{conditionally-strong quasi-stationary time}  (CSQST ) 
if for any $y \in A$, and $t\ge 0$
\be{defCSQST'}
 \P(X^{\alpha}_t=y,\, \t^{\alpha}_*=t)=\m^*(y)\P(\t^{\alpha}_*=t<\t^{\alpha}_G).
 \ee
or, in other words, for any $y\in A$ and $t\ge 0$
\be{CSQST}
\P \left( X^\alpha_t = y, \tau^\alpha_*=t \ | \ t < \tau^\alpha_G \right) 
= 
\mu^*(y) \P \left(  \tau^\alpha_*=t \ | \ t < \tau^\alpha_G \right) 
\ee
\end{itemize}
 
 %%%%%%%%
 \subsubsection{The local chain $\widetilde X_t$ on $A$}
\label{localch}
In this subsection we construct an ergodic Markov chain $\widetilde X_t$ on $A$, that we call 
the {\it local chain}. 

Many  dynamics have been used in the literature to describe the local behavior of the process $X_t$ on $A$.
Examples are the reflected process or the conditioned process (see for instance \cite{FMNSS}, \cite{BG}).

We use here a  local chain  $\widetilde X_t$ constructed by means of the 
 right eigenvector  of $[P]_A$ corresponding to $\l$.
 This construction is related to the Doob h-transform of $[P]_A$ (see for instance \cite{LPW}).
This chain $\widetilde X_t$ is also related to the  ``reversed chain" in {Darroch-Seneta}, introduced
in \cite{DS} while considering
the large time asymptotics. 

\bigskip

The construction is the following:
 by the Perron-Frobenius theorem 
  there exists a unique positive right eigenvector $\g$ of $[P]_A$ corresponding to $\l$, i.e.,
  \be{eqmqs1}
  [P]_A\g=\l \g\qquad \hbox{ with normalization }\qquad \m^*\cdot\g=1.
  \ee
This eigenvector is related to the asymptotic ratios
of the survival probabilities {(see eg \cite{CMS})}
$$
\lim_{t\to \infty}\frac{\P(\t^x_G>t)}{\P(\t^y_G>t)}=\frac{\g(x)}{\g(y)}\qquad x,y\in A.
$$

For any $x,y\in A$, define the stochastic matrix
\be{tildeP}
\widetilde P(x,y):=\frac{\g(y)}{\g(x)} \frac{P(x,y)}{\l}.
\ee
Notice that $\widetilde P$ is a primitive matrix.
Let $\nu$ be its invariant measure
$$
\sum_{x\in A}\nu(x)\widetilde P(x,y)=\nu(y)=\sum_{x\in A}\nu(x) \frac{\g(y)}{\g(x)} \frac{P(x,y)}{\l}
$$
it is easy to see that 
$$\g(x)=\frac{\nu(x)}{\m^*(x)},\qquad \forall x\in A$$
For the chain $\widetilde X_t$ we define
$$\tilde s^x(t,y):=1-\frac{\widetilde P^t(x,y)}{\nu(y)}$$
$$\tilde s^x(t)= sep(\tilde \m^x_t,\nu)=\sup_{y\in A}\tilde s^x(t,y),\quad\tilde s(t):=\sup_{x\in A}\tilde s^x(t).$$
Note that $\tilde s^x(t)\in[0,1]$. Moreover, since $\widetilde P$ is a primitive matrix, it is well known (see for instance \cite{AD1}, Lemma 3.7) that
 $\tilde s(t)$ has 
 the sub-multiplicative property:
$$\tilde s(t+u)\le \tilde s(t)\tilde s(u).$$
This implies in particular an exponential decay in time of $\tilde s(t)$.

\bigskip

The relation between the local chain and the original chain $X_t$ on $\cX$
is given by the definition (\ref{tildeP}) and more generally by
\be{tildePt}
\widetilde P^t(x,y)=\frac{\g(y)}{\g(x)} \frac{P^t(x,y)}{\l^t}\qquad\forall t\ge 0.
\ee
%%%%%%%%%%
\subsubsection{Preliminary remarks}
\label{prel_rem}
We can use  this last relation 
to obtain a {\it rough estimate about the absorption time} $\t_G$.
We give here this simple calculation in order to point out the dependence on the initial distribution $\a$
of the distribution of $\t^{\a}_G$ by means 
 of  a  {\it time shift} defined by
 \be{da}
\d_{\a}:=\log_\l\big(\a\cdot\g\big) 
\ee

 We will show that it is natural to associate to every initial measure $\a$
the following measure $\tilde\a$ for the local chain $\widetilde X_t$:
$$
\tilde \a(x) := \frac{\a(x)\g(x)}{\a\cdot\g}.
$$

Indeed,
 $$
\P(\t^{\a}_G>t)=\sum_{y\in A}\sum_{x\in A}\a(x) {P^t(x,y)}=$$
$$
\sum_{y\in A}\sum_{x\in A}\a(x) {\g(x)\l^t\m^*(y)\frac{\widetilde P^t(x,y)}{\nu(y)}}=
$$
\be{ub_rozzo}
\l^t\sum_{x\in A}\a(x) \g(x)\sum_{y\in A}\m^*(y)(1-\tilde s^x(t,y))=$$
$$\l^{t+\d_\a}\Big(1-\sum_{y\in A}\m^*(y)\tilde s^{\tilde\a}(t,y)\Big)\ge
\l^{t+\d_\a}\Big(1-\tilde s^{\tilde\a}(t)\Big)
\ee
with
\be{stildeatilde}
\tilde s^{\tilde\a}(t,y) := \sum_{x\in A}\tilde\alpha(x)\tilde s^x(t,y)\qquad \hbox{ and} \qquad \tilde s^{\tilde\a}(t):=\sup_{y\in A}\tilde s^{\tilde\a}(t,y)
\ee
Note that from (\ref{ub_rozzo}) we obtain for any initial distribution $\a$
\be{notada}
\l^{t+\d_\a}(1-\tilde s^{\tilde\a}(t))\le 1\qquad \forall t\ge 0.
\ee
To obtain an upper bound on $\P(\t^{\a}_G>t)$, 
we can consider the minimal  strong stationary time $\tilde \t^x_\nu$ such that
$$
\P(\widetilde X^x_t=y,\; \tilde \t^x_\nu=t)=\nu(y) \P(\tilde \t^x_\nu=t)
$$
 with
$$
\P(\tilde\t^x_\nu>t)=\tilde s^x(t).
$$
Similarly, we have
$$
\P(\t^{\a}_G>t)=\sum_{y\in A}\sum_{x\in A}\a(x) \g(x)\l^t\frac{\m^*(y)}{\nu(y)}\P(\widetilde X^x_t=y, \tilde\t^x_\nu\le t)+
$$
$$
\sum_{y\in A}\sum_{x\in A}\a(x) \g(x)\l^t\frac{\m^*(y)}{\nu(y)}\P(\widetilde X^x_t=y, \tilde\t^x_\nu> t)\le
$$
$$
\sum_{y\in A}\sum_{x\in A}\a(x) \g(x)\l^t\frac{\m^*(y)}{\nu(y)}\nu(y)\P(\tilde\t^x_\nu\le t)+
\frac{1}{\min_y\g(y)}\sum_{x\in A}\a(x) \g(x)\l^t\P(\tilde\t^x_\nu>t)
$$
$$
=\l^{t+\d_{\a}}\Big[1+\tilde s^{\tilde\a}(t)\big(\frac{1}{\min_y\g(y)} -1\big) \Big].
$$
This quantity could be much larger that $1$, since $\frac{1}{\min_y\g(y)}\ge 1$,
and so this estimate from above on the distribution of $\t^\a_G$ is quite rough.  However, we have to note
that this factor is independent of time so that, for large $t$,  due to the exponential decay of $\tilde s(t)$,
and so
of $\tilde s^{\tilde\a}(t)$,  the estimate is not trivial. 

It is interesting to notice that for sufficiently large $t$, say $t\ge \min\{n\in \mathbb N:\; n+\d_\a\ge 0\}$, the separation $\tilde s^{\tilde\a}(t)$ has a straightforward meaning 
for the $X_t$ process:
it is related to the separation between the measure $\mu^\alpha_t$ and the evolution starting from the quasi-stationary measure corrected with a time-shift $\delta_{\alpha}$, namely the measure 
$${\m^{\m^*}_{t+\d_{\a}}(y)}:=\left\{\begin{array}{c}
\lambda^{t+\delta_{\alpha}}\mu^*(y) \ \ \ \ \ \text{ if } y \in A\\
1 - \lambda^{t+\delta_{\alpha}} \omega(y) \ \text{ if } y \in G
\end{array}\right. ,$$ 
where $\omega$ is the hitting distribution defined in \eqref{defomega}. More precisely,
by the definition of the process $\widetilde X_t$, we have for any $y\in A$:
\be{s=stilde}
\tilde s^{\tilde\a}(t,y):= 1 - \sum_{x\in A}\tilde \a(x)\frac{\tilde P^t(x,y)}{\nu(y)}=1-\frac{\m^{\a}_t(y)}{\m^{\m^*}_{t+\d_{\a}}(y)}.
\ee
This means that

$$
\sum_{x\in A}\tilde \a(x)\frac{\widetilde P^t(x,y)}{\nu(y)}=\frac{\m^\a_t(y)}{\m^{\m^*}_{t+\d_\a}(y)},\qquad y\in A
$$
 so that the convergence to equilibrium of the local chain $\widetilde X^{\tilde\a}_t$ controls the convergence of the 
chain $X^\a_t$ to the evolution starting from the quasi-stationary measure corrected with a time-shift, $\m^{\m^*}_{t+\d_\a}$,
as far as its permanence in the set $A$ is concerned.
This is the reason why the local chain $\widetilde X_t$ is crucial in our discussion.

%%%%%%%%%%%%
\bigskip

%%%%%%%%%%%%%%
%%%%%%%%%%%%%% 
 \subsection{Main results}
\label{mainres} 
 
We collect in this section our main results on conditionally strong quasi stationary
 times (CSQST) and their application to control the distribution of the hitting time $\t_G$.
 
 From the  definition of  CSQST  (\ref{defCSQST'}) we can prove the following:
\bp{p1}
For any initial distribution $\a$ on $A$ and
 for any    conditionally strong quasi stationary
 time  $\t^{\a}_*$   we have for any $t\ge 0$:
 $$
 \P(\t^{\a}_*\le t<\t^{\a}_G)=\sum_{u\le t}\l^{t-u}\P(\t^{\a}_*=u<\t^{\a}_G)\le \l^{t+\d_{\a}}(1-\tilde s^{\tilde\a}(t)).
 $$

\ep

Proposition \ref{p1}  suggests a new notion of minimality.
\bd{d_min_CSQST}
For any initial distribution $\a$ on $A$ a
    conditionally strong quasi stationary
 time  $\t^{\a}_*$   is \emph{minimal } if for any $t\ge 0$:
 $$
 \P(\t^{\a}_*\le t<\t^{\a}_G)= \l^{t+\d_{\a}}(1-\tilde s^{\tilde\a}(t)).
 $$
\ed

 The existence of minimal conditionally strong quasi-stationary times is given by the following Theorem.
 %%%%
 \bt{main1}
  For any initial distribution $\a$ on $A$,   there exists a \emph{minimal conditionally strong quasi stationary
 time} $\t^{\a}_*$  
 such that for any $t>0$
$$\P(\t^{\a}_*=t<\t^{\a}_G)= \l^{t+\d_{\a}}(\tilde s^{\tilde\a}(t-1)-\tilde s^{\tilde\a}(t)).$$
 \et

 Note that in particular for a minimal conditionally strong quasi stationary
 time we have for $ t\ge 0$
$$\P(\t^{\a}_*>t, \; \t^{\a}_*<\t^{\a}_G)= \sum_{u>t}\l^{u+\d_{\a}}(\tilde s^{\tilde\a}(u-1)-\tilde s^{\tilde\a}(u))
\le \l^{t+\d_\a}\tilde s^{\tilde\a}(t).$$
\bigskip

Let $\t^\a_*$ be a minimal CSQST and define  
$$\t^{\a}_{*,G}={\t^{\a}_G\wedge \t^{\a}_*}.$$
This time plays the role of  {\it local relaxation time} or {\it metastability time}, like $\t_{m\cup G}$  in the metastable hypothesis (\ref{rar0}).
It is a sub-multiplicative time:
  \bt{submult}
  If $\tau^\alpha_*$ is a minimal CSQST, then for any positive $u$ and $v$ 
  \be{eq:submult}
    \sup_{\alpha} \P\Big( \t^{\a}_{*,G}> u+v \Big) \le
    \sup_{\alpha} \P\Big( \t^{\a}_{*,G}> u \Big) 
    \sup_{\alpha} \P\Big( \t^{\a}_{*,G}> v \Big) .
  \ee
  \et

The local relaxation time $\t^{\a}_{*,G}$ is a key ingredient in the following representation formula:
 \bt{main2}
 For any initial distribution $\a$ on $A$,  if  $\t^{\a}_*$   is a minimal conditionally strong quasi stationary
 time,  we have, for any $t\ge 0$
 \be{main2_1}
  \P\Big(\t^{\a}_G>t\Big)=\l^{t+\d_{\a}}(1-\tilde s^{\tilde\a}(t))+\P\Big( \t^{\a}_{*,G}> t\Big)
 \ee

 Moreover, for any $y\in G$, we have
 \be{uscita}
 \P\Big(X^{\a}_{\t^{\a}_G}=y\Big)=\P\Big( \t^{\a}_G< \t^{\a}_*, \; X^{\a}_{\t^{\a}_G}=y\Big)+
 \omega(y)\P\Big( \t^{\a}_G> \t^{\a}_*\Big),
\ee
where $\omega$
is the hitting distribution starting from $\m^*$ (see equation (\ref{defomega})).
 \et 
  %%%%
  \bigskip
  
  This theorem provides a    control on the convergence  to an exponential distribution
  for the hitting time $\t_G$ and  on the hitting distribution
  and it gives a probabilistic interpretation of the errors in the exponential approximation of 
  $ \P\Big(\t^{\a}_G>t\Big)$ in terms of conditionally-strong quasi-stationary times.

  In order to obtain a multiplicative bound on the exponential distribution, it is useful to rewrite eq. \eqref{main2_1} as
  \be{m1}
     \frac{\P\Big(\t^{\a}_G>t\Big)}{\lambda^{t+\delta_{\alpha}}} -1 =
     - \tilde s^{\tilde\a}(t)+\lambda^{-t-\delta_\alpha}\P\Big( \t^{\a}_{*,G}> t\Big).
  \ee
  The first error term $- \tilde s^{\tilde\a}(t)$ decays exponentially fast in $t$
  and it will be easy to deal with;
  the second error term $\lambda^{-t-\delta_\alpha}\P\Big( \t^{\a}_{*,G}> t\Big)$ will decay faster than the 
  leading term only under suitable metastability hypotheses.
  
  Let $\tau^\alpha_*$ be a minimal CSQST. Define the {\it mean metastability time}:
  \be{mmt}
  R:=\sup_{\alpha} \E  \Big(\tau^\alpha_{*,G}\Big)
  \ee
  and the {\it mean relaxation time}:
  \be{mrt}
  T:=(1-\lambda)^{-1}=\E \Big( \tau^{\m^*}_G\Big)
  \ee
  
 \bd{dmmh}
 We call {\emph{ mean metastability hypothesis}} with rate $a$ the condition
  \be{mmh}
     \frac{R}{T} =\l^T e^{-a},
  \ee
  for some $a>0$.
  \ed
  Note that for   $
  e^{-1/\l}\le \l^T\le e^{-1}
  $
  that are strict bounds if $T$ is large, as in metastable situations.

  \bt{mainexp}
     Under the mean metastability hypothesis given in Definition \ref{dmmh}, for any initial measure $\a$ on $A$
      and 
    for any positive integer $n$,
    \be{eq:mainexp}
      \left| \frac{\P\Big( \t^{\a}_{G}> n T \Big)}{\lambda^{nT+\delta_\alpha} } -1\right| <{e^{-an}}\l^{-\d_\a}
      \frac{e^{1/\l}}{1-e^{-a}}.
    \ee
  \et
 %%%%%
 Note that the mean metastability hypothesis (\ref{mmh}) does not exclude the existence of starting states  $x\in A$ 
 from which the process reaches $G$ in a very short time with high probability. When the starting distribution
 $\a$ is concentrated on such states, we expect to have $\l^{\d_\a}$ very small. 
 This implies that  Theorem \ref{mainexp} provides a sharp result, in the case of small $n$, only if
the parameter $a$ in (\ref{mmh}) is sufficiently large and $\l^{\d_\a}$ is not too small. More precisely, if the
 support of the starting measure $\a$ is contained in a ``basin of attraction of the metastable state" defined  for instance
 as
 (see \cite{FMNSS})
 $$
 B:=\Big\{x\in A:\; \P(\t^x_G>2R)>3/4\Big\},
 $$
 we can give a very rough estimate  $\l^{\d_\a}\ge 1/4$.
  
Indeed by using the trivial estimate
 $$
 \P(\t^\a_G>nT)=\sum_{y\in A}\P(X^\a_{2R}=y)\P(\t^y_G>nT-2R)\ge  \l^{nT-2R}\P(\t^\a_*\le 2R<\t^\a_G)=$$
 $$
 \l^{nT-2R}\Big[\P(\t^\a_G>2R)-\P(\t^\a_{*,G}>2R)\Big],
 $$
 by
 Theorem \ref{mainexp} and the Markov inequality we get
$$\l^{\d_\a}= \lim_{n\to\infty}\frac{\P(\t^\a_G>nT)}{\l^{nT}}\ge \l^{-2R}\Big[\P(\t^\a_G>2R)-\P(\t^\a_{*,G}>2R)\Big]
 \ge 1/4.$$
 
 %%%%%
  In many applications it is interesting to study the behavior of the process on an intermediate time-scale $S$, say  $R\ll S\ll T$.  
  The process has an {\emph{early exponential behavior}} if equation (\ref{eq:mainexp}) holds by replacing $T$ with $S$.
  In \cite{FMNS} the early exponential behavior of the first hitting time is proved in a  particular case, with
  a particular starting configuration. 
  In our setting, we can study the early behavior starting from a general measure $\a$ under the mean metastability hypothesis with very large rate $a$.

\bigskip

%%%%%%%%%%%%%
%%%%%%%%%%%%%%%
 %%%%%%%%%%%%%%%%
 %%%%%%%%%%%%%%%%%
%%%%%
%%%%%%%
\section{Proofs}
\label{proofth}
\subsection{Tracking process} \label{trp}
%\color{blue}
In order to prove the existence of a minimal CSQST, we introduce
 an auxiliary \emph{tracking process}. 
 The construction   is inspired to  \cite{AD1} and 
\cite{DF90}, where  the existence of strong stationary times is proved.
The idea is to duplicate the state space into two  layers and to define a process on this larger state space with
a jump probability from one layer to the other one. In order to have a general construction, we introduce first a control 
function to define the jump rate.
\bd{defmt}
Let $\Z_{\ge -1}$ denotes  integers larger or equal to $-1$.
The function   $m(t): \Z_{\ge -1} \longrightarrow [0,1]$ is a {\emph{ control function for the process starting at $\a$}} if it is a monotonic decreasing function with 
$$
m(t)\ge \tilde s^{\tilde\a}(t)  \text{ for } t \ge 0, \qquad m(-1)=1.
$$
\ed

Given a control function $m(t)$  for every $z\in A$, we define the following {\it jump probabilities} for any $t\ge 0$ 
\begin{equation}
J^{\a}(t,z):=\frac{m(t-1)-m(t)}{m(t-1)-\tilde s^{\tilde\a}(t,z)},\label{eq:defJ}
\end{equation}
with the convention $0/0=0$.
Since $m(t-1)\ge m(t)\ge s^{\tilde\a}(t,z)$, we have $J^{\a}(t,z)\in[0,1]$
for any $ z\in A$ and any $t$.
For any $t\ge 0$ and any $z\in G$ we define
$$
J^{\a}(t,z)\equiv J^{\a}(t,G):= 1.
$$
\bd{defY}
On the state space $\boldsymbol {\mathcal{ X}}:=\mathcal{X}\times\{0,1\}$, consider the transition matrix 
\begin{eqnarray}\label{defQ}
&&Q^{\a}_t \left( (y,0),(z,0)\right):=P(y,z)\Big(1-J^{\a}(t,z)\Big), \nonumber\\  
&&Q^{\a}_t\left({(y,0),(z,1)}\right):=P(y,z)J^{\a}(t,z),\nonumber\\ 
&&Q^{\a}_t \left({(y,1),(z,e)}\right):=P(y,z){\mathds 1}_{\{e=1\}},
\end{eqnarray}
where $e\in \{0;1\}$;
also, consider the initial distribution ${{\boldsymbol\alpha}}$  on the two layers of $\boldsymbol {\mathcal{ X}}$, defined,
for any $x\in A$, as
\begin{eqnarray}\label{inizY1}
{{\boldsymbol\alpha}}(x,0)&:=&\a(x)\big(1-J^\a(0,x)\big)=\a(x)-\l^{\d_\a}\m^*(x) (1-m(0)), \nonumber\\
{{\boldsymbol\alpha}}(x,1)&:=&\a(x)J^\a(0,x)=\l^{\d_\a}\m^*(x) (1-m(0)).
\end{eqnarray}

We define the {\emph{ tracking process}} ${\boldsymbol X}^{{{\boldsymbol\alpha}}}_t $ via 
$$
\P \left(\bigcap_{u=0}^t ({\boldsymbol X}^{{\boldsymbol\alpha}}_u = \boldsymbol{y}_u) \right) = {{\boldsymbol\alpha}}(\boldsymbol y_0) \prod_{u=0}^{t-1} Q^\alpha_u (\boldsymbol y_u,\boldsymbol y_{u+1}) 
$$
with $\boldsymbol y_u\in \boldsymbol {\mathcal{ X}}$ for any $u\le t$.
\ed
By \eqref{defQ}, \eqref{inizY1} it is immediate to see 
 that the marginal distribution of ${\boldsymbol X}_t^{{{\boldsymbol\alpha}}}$ on $\cX$ corresponds to the distribution
of $X^{\a}_t$, so that we can study each event defined for the process $X^{\a}_t$ in terms
of sets of paths of the process ${\boldsymbol X}^{{{\boldsymbol\alpha}}}_t$. 
For this reason, with an  abuse of notation, we denote with the same symbol
$\P$ the probability of events defined in terms of  the process  ${\boldsymbol X}^{{{\boldsymbol\alpha}}}_t$.

%%%%
%%%%% 
Notice that unlike the process defined in \cite{AD1} and 
\cite{DF90} for the strong stationary times, 
 in our definition  of the jump rates we use
 the separation $\tilde s^{\tilde\a}$ for the process ${\widetilde X^{\tilde\a}}_t$, defined
on $\widetilde \cX$.  

%%%%
We want also to note that the starting measure $\alpha$ appears as a parameter in the definition of the transition matrix $Q$ (see \eqref{eq:defJ}, \eqref{defQ}), 
 the process is time-inhomogeneous  and Markov property does not hold. 
However, we can get rid of this dependence and recover a sort of semigroup property by considering a suitable conditioning of the process ${\boldsymbol X}^{{\boldsymbol\alpha}}_t$. We will clarify this point, that represents a crucial ingredient in our approach, in Section \ref{ephemeral}.

We will be interested to the process ${\boldsymbol X}^{{\boldsymbol\alpha}}_t$  up to its first hitting to the
set $\cX \times \{1\}$, i.e. for $t\le \tau^{{{\boldsymbol\alpha}}}_1$ with
\be{tau1}
\tau^{{{\boldsymbol\alpha}}}_1:=
\tau^{{{\boldsymbol\alpha}}}_{\cX \times \{1\}}=\min\{t\ge0\:;\: {\boldsymbol X}_{t}^{{{\boldsymbol\alpha}}}=(y,1)\:\text{for some }y\in\mathcal{X}\},
\ee
indeed, we  prove 
that  $\tau^{{{\boldsymbol\alpha}}}_1$ is a conditionally strong quasi-stationary time.

This construction of a CSQST is quite implicit, for it requires the knowledge of the
separation $\tilde s^{\tilde \a}(t)$ at any time, which in general is very hard to obtain.
In this paper we use CSQST as a theoretical tool and we are not concerned with their explicit construction.
However, it is well-known that in some systems the separation can be estimated with the
distribution of a hitting time to a suitable halting state (see e.g. \cite{LPW}).
In the example in section \ref{ex}, we exploit this idea to construct explicitly 
a (non-minimal) CSQST.

\color{black}
%%%%%%%%%%%%%%%%%%%%%%%%%%%%%%%%%%%%%%%%%%%%%%%%%%%%%%%%%%%%%%%%%%%%%

\color{black}

\subsection{ Conditionally strong quasi stationary times (CSQST)}
\label{s_main1}

In this section we prove Proposition \ref{p1} and Theorem \ref{main1}.

Let us start by proving Proposition \ref{p1}:
by the definition of CSQST we have for any $y\in A$
\be{p1'1}
 \P(X^{\a}_t=y,\, \t^{\a}_*=u)=\m^*(y)\P(\t^{\a}_*=u<\t^{\a}_G),\qquad \hbox{for any }\; u\ge 0 
 \ee

If (\ref{p1'1}) hods for any $u\in [0,t]$ then 
we have:
$$
 \P(X^{\a}_t=y,\, \t^{\a}_*\le t)= \sum_{u\le t}\sum_{z\in A}\P(X^{\a}_u=z,\, \t^{\a}_*=u)P^{t-u}(z,y) =
$$
\be{defCSQSTle}
\sum_{u\le t}\sum_{z\in A}\m^*(z)\P(\t^{\a}_*=u<\t^{\a}_G)P^{t-u}(z,y)=\m^*(y)\sum_{u\le t}\l^{t-u}\P(\t^{\a}_*=u<\t^{\a}_G)
 \ee
 and by summing over $y\in A$:
\be{t*<t}
\P(\t^{\a}_*\le t <\t^{\a}_G)=\sum_{u\le t}\l^{t-u}\P(\t^{\a}_*=u<\t^{\a}_G).
\ee
%%%%%%

Moreover for any $y\in A$ we have
\begin{eqnarray*}
\m^{\a}_t(y)
&\ge& 
\P(\t^{\a}_*\le t, X^{\a}_t=y)=
 \sum_{u\le t}\sum_{z\in A}\P(\t^{\a}_*=u, X^{\a}_u=z)P^{t-u}(z,y)=\\
&&\l^t
 \sum_{u\le t}\l^{-u}\P(\t^{\a}_*=u<\t^{\a}_G)\m^*(y)
\end{eqnarray*}
so that by (\ref{s=stilde})
$$
\frac{\m^{\a}_t(y)}{\l^t\m^*(y)}=\l^{\d_{\a}}(1-\tilde s^{\tilde\a}(t,y))\ge  \sum_{u\le t}\l^{-u}\P(\t^{\a}_*=u<\t^{\a}_G)=
\P(\t^{\a}_*\le t<\t^{\a}_G)\l^{-t}
$$
since this holds for any $y\in A$  we get
 \be{bo1}
 \P(\t^{\a}_*\le t<\t^{\a}_G)=\sum_{u\le t}\l^{t-u}\P(\t^{\a}_*=u<\t^{\a}_G)
 \le \l^{t+\d_{\a}}(1-\tilde s^{\tilde\a}(t)).
 \ee
\cvd
\bigskip

%%%%%
We prove a stronger version of Theorem \ref{main1}.
%%%%%
Indeed by choosing the control function $m(t)=\tilde s^{\tilde\a}(t)$, Theorem \ref{main1}
immediately follows by:
\bt{main1'}
For any initial distribution $\a$ on $A$
  and for any control function $m(t)$ 
  there exists a \emph{ conditionally strong quasi stationary
 time} $\t^{\a}_*$  such that
 $$
 \P(\t^{\a}_*\le t<\t^{\a}_G)= \l^{t+\d_{\a}}(1-m(t))\qquad \hbox{ for all }\; t\ge 0
 $$
with 
$$\P(\t^{\a}_*=t<\t^{\a}_G)= \l^{t+\d_{\a}}(m(t-1)-m(t)).$$
 \et

%%%%
To prove  Theorem \ref{main1'} consider now  the  tracking process
defined in Section \ref{trp} 
and  the hitting time
\[
\tau^{{\boldsymbol\alpha}}_1:=
\tau^{{\boldsymbol\alpha}}_{\cX \times \{1\}}=\min\{t\ge0\:;\: {\boldsymbol X}_{t}^{{\boldsymbol\alpha}}=(y,1)\:\text{for some }y\in\mathcal{X}\}.
\]

\color{black}

We will prove that  $\tau^{{\boldsymbol\alpha}}_1$ satisfies the
following condition for any $t\ge 0$:
$$
{\mathcal C}(t):=
\begin{cases}
 \P(X^{\a}_t=y,\, \t^{{\boldsymbol\alpha}}_1=t)=\m^*(y)\P(\t^{{\boldsymbol\alpha}}_1=t<\t^\a_G)\quad \hbox{ for any } y\in A\\
\P(\tau^{{\boldsymbol\alpha}}_1=t<\t^\a_G)=\l^{t+\d_\a}\big(m(t-1)-m(t)\big)
\end{cases}
$$
If  ${\mathcal C}(u)$ is verified for any $u\le t$, we can conclude 
$$
\P(\tau^{{\boldsymbol\alpha}}_1\le t<\t^\a_G)=\sum_{u\in[0, t]}\l^{t-u}\P(\tau^{{\boldsymbol\alpha}}_1=u<\t^\a_G)= \l^{t+\d_\a}\big(1-m(t)\big).
$$
\bigskip

In order to prove that $\tau^{{\boldsymbol\alpha}}_1$ satisfies ${\mathcal C}(t)$  for all $t\ge 0$ we proceed by induction on $t$.
For $t=0$,  by the definition of the initial distribution ${\boldsymbol\alpha}$ in definition
\ref{defY} we immediately verify
 ${\mathcal C}(0)$. 
 Indeed, for $y\in A$ we get
 $$\P(X^{\a}_{0}=y,\, \t^{{\boldsymbol\alpha}}_1={0})=
\P({\boldsymbol X}^{{\boldsymbol\alpha}}_{0}=(y,1))=\l^{\d_\a}\m^*(y)\big(1-m({0})\big)=$$
$$\m^*(y) \P(\t^{{\boldsymbol\alpha}}_1={0}<\t^\a_G)$$
with
$$
\P(\t^{{\boldsymbol\alpha}}_1={0}<\t^\a_G)=\sum_{x\in A}\P({\boldsymbol X}^{{\boldsymbol\alpha}}_{0}=(x,1))=\l^{\d_\a}\big(1-m({0})\big).
$$

To prove the induction step we use the following:
\bl{l001}
If for any $u\in[0,t]$  and for any $y\in A$ we have
$$
\P\big( X^{\a}_u=y,\;\t^{{\boldsymbol\alpha}}_1=u\big)={\m}^*(y)\P(\t^{{\boldsymbol\alpha}}_1=u<\t^\a_G)
$$
then, for any $z\in A$,
$$
\P\big( {\boldsymbol X}^{{\boldsymbol\alpha}}_t=(z,1)\big)={\m}^*(z)\P(\t^{{\boldsymbol\alpha}}_1\le t<\t^\a_G).
$$
\el
\bpr
Note first that under the hypothesis of the Lemma, by (\ref{defCSQSTle}) we get
\be{lemma1}
\P(\t^{{\boldsymbol\alpha}}_1\le t<\t^\a_G)=\sum_{u\le t}\l^{t-u}\P(\t^{{\boldsymbol\alpha}}_1=u<\t^\a_G).
\ee
We have
$$
\P\big( {\boldsymbol X}^{{\boldsymbol\alpha}}_t=(z,1)\big)=\sum_{u\le t}\sum_{y\in A}\P\big( {\boldsymbol X}^{{\boldsymbol\alpha}}_t=(z,1),\; \t^{{\boldsymbol\alpha}}_1=u,\: X^{\a}_u=y\big)=
$$
$$\sum_{u\le t}\sum_{y\in A}\P(\t^{{\boldsymbol\alpha}}_1= u,\; X^\a_u=y)P^{t-u}(y,z)
=\sum_{u\le t}\P(\t^{{\boldsymbol\alpha}}_1= u)\sum_{y\in A}\m^*(y)P^{t-u}(y,z)=
$$
$$
=\m^*(z)\sum_{u\le t}\l^{t-u}\P(\t^{{\boldsymbol\alpha}}_1=u<\t^\a_G)=\m^*(z)\P(\t^{{\boldsymbol\alpha}}_1\le t<\t^\a_G).
$$
\epr

Suppose now that ${\mathcal C}(u)$  holds for $u\in[0, t]$.
By using then Lemma  \ref{l001} we get
$$
\P\big( X^{\a}_{t+1}=y, \t^{{\boldsymbol\alpha}}_1={t+1} \big)=\sum_{z\in\cX}\P\big({\boldsymbol X}^{{\boldsymbol\alpha}}_{t+1}=(y,1)|{\boldsymbol X}^{{\boldsymbol\alpha}}_{t}=(z,0)\big)\P\big({\boldsymbol X}^{{\boldsymbol\alpha}}_{t}=(z,0)\big)=
$$
$$
\sum_{z\in\cX}P(z,y)J^{\a}({t+1},y)\big[ \m^{\a}_{t}(z)-\P({\boldsymbol X}^{{\boldsymbol\alpha}}_{t}=(z,1)) \big]=$$
\be{passoind1}
J^{\a}({t+1},y)\big[\m^{\a}_{t+1}(y)-\sum_{z\in\cX}{\m}^*(z)P(z,y)\P(\t^{{\boldsymbol\alpha}}_1\le t<\t^\a_G)\big].
\ee
Since ${\mathcal C}(u)$  holds for $u\in[0, t]$ we have
$$
\P(\t^{{\boldsymbol\alpha}}_1\le t<\t^\a_G)=\sum_{u\le t}\l^{t-u}\P(\t^{{\boldsymbol\alpha}}_1=u<\t^\a_G)= \l^{t+\d_\a}(1-m(t))
$$
Recalling that, by (\ref{s=stilde}),
$$
\m^{\a}_{t+1}(y)=\l^{t+1+\d_\a}\m^*(y)[1-\tilde s^{\tilde\a}({t+1},y)],
$$
we obtain
$$
\P\big( X^{\a}_{t+1}=y, \t^{{\boldsymbol\alpha}}_1={t+1} \big)=J^{\a}({t+1},y)
\l^{t+1+\d_\a}\m^*(y)\big[1-\tilde s^{\tilde\a}({t+1},y)-(1-m(t))\big].
$$
By using the definition of $J^\a({t+1},y)$, we get
$$
\P\big( X^{\a}_{t+1}=y, \t^{{\boldsymbol\alpha}}_1={t+1} \big)=[m(t)-m(t+1)]\l^{t+1+\d_\a}\m^*(y)
$$
so that, by summing on $y\in A$
$$
\P\big( \t^{{\boldsymbol\alpha}}_1={t+1}<\t^\a_G \big)=[m(t)-m(t+1)]\l^{t+1+\d_\a},
$$
we show that ${\mathcal C}(t+1)$ holds, concluding the proof of Theorem \ref{main1'} and so of Theorem \ref{main1}.\cvd

\label{tpem}
%%%%%%%%%
%%%%%%%%%

\color{black}
%%%%%%%%%
%%%%%%%%%
\subsection{Ephemeral measure}\label{ephemeral}

In this section we describe the behavior of the process before $\tau_{1}^{{{\boldsymbol\alpha}}}$.
We call this behavior  ``ephemeral" since in
 metastable situations $\tau_{1}^{{{\boldsymbol\alpha}}}$ is typically much smaller than $\t_G$. 
 
 Consider  the tracking process before $\tau_{1}^{{{\boldsymbol\alpha}}}$,
more precisely,  the conditioned measure on $\cX\times \{0\}$ obtained by the process ${\boldsymbol X}_{t}^{{{\boldsymbol\alpha}}}$
conditioned to the layer $\{0\}$:
\bd{deph}
The measure
\be{def_eph}
\Phi_{t}^{{\boldsymbol\alpha}}(x):=\mathbb{P}\left({\boldsymbol X}_{t}^{{\boldsymbol\alpha}}=(x,0)\:|\vphantom{\:}\tau_{1}^{{{\boldsymbol\alpha}}}>t\right).
\ee
is called  the {\emph ephemeral measure}. 
 \ed
With a slight abuse of notation, we consider this ephemeral measure either as a measure on 
 $\boldsymbol {\mathcal{ X}}$ (with support in $A\times \{0\}$) or as a measure in $A$.
 Recalling that ${{\boldsymbol\alpha}}(x,0)=\a(x)\Big(1-J^\a(0,x)\Big)$ we get for the ephemeral measure:
 $$
 \Phi_{t}^{{\boldsymbol\alpha}}(x)= \frac{1}{\mathbb{P}\left(\tau_{1}^{{\boldsymbol\alpha}}>t\right)}
 \sum_{x_0\in A}\a(x_0) \Big(1-J^\a(0,x_0)\Big)\times$$
 $$\sum_{x_1,...,x_{t-1}}\Big[\prod_{s=1}^{t-1}
P(x_{s-1},x_s)\Big(1-J^\a(s,x_s)\Big)\Big]
P(x_{t-1},y)\Big(1-J^\a(t,y)\Big)
$$
We will prove the following ``Markov-like" properties for the tracking process and for the hitting
time $\tau_{1}^{{{\boldsymbol\alpha}}}$:
\bp{p_eph}
Consider the tracking process ${\boldsymbol X}^{{\boldsymbol\alpha}}_t$ starting at $\boldsymbol\a$ and with control function 
\be{mt}
m(t)=
 \tilde s^{\tilde \alpha}(t) \quad\forall t \ge 0
\ee
then for any $x\in A$ 
$$
\Phi_{t+u}^{{\boldsymbol\alpha}}(x)=\Phi_u^{\Phi_{t}^{{\boldsymbol\alpha}}}(x)
$$
\ep

\bp{p_tau1}
Consider the tracking process ${\boldsymbol X}^{{\boldsymbol\alpha}}_t$  with control function $m(t)=\tilde s^{\tilde \alpha}(t), \quad\forall t \ge 0$,
then 
$$
\mathbb{P}(\tau_{1}^{{\boldsymbol\alpha}}>t+u) =  \mathbb{P}(\tau_{1}^{{\boldsymbol\alpha}}>t)\mathbb{P}(\tau_{1}^{\Phi_{t}^{\boldsymbol\alpha}}>u)
$$
\ep
From this Proposition the submultiplicativity property of Theorem \ref{submult} easily follows (see Section \ref{sub}).

To prove these propositions, we introduce two technical lemmas to obtain the crucial property
on the jump rates given in Lemma \ref{Jf}.

Recalling from section \ref{localch} that 
\[
\tilde{\alpha}(x)=\frac{\alpha(x)\gamma(x)}{\alpha\cdot\gamma},
\]
 we denote
by $\tilde{\Phi}_{t}^{\tilde\alpha}$ the measure
\[
\tilde{\Phi}_{t}^{\tilde\alpha}(x)=\frac{\Phi_{t}^{{\boldsymbol\alpha}}(x)\gamma(x)}{\Phi_{t}^{{\boldsymbol\alpha}}\cdot\gamma}.
\]

\bl{KH}

There exist two functions $K$ and $H$, that depend on $\alpha$
and on $t$ but not on $x$, such that, for any $y\in A$, 

\[
\tilde{\Phi}_{t}^{\tilde{\alpha}}(y)=K\sum_{x\in A}\tilde\a(x)\tilde P^t(x,y)-H\nu(y).
\]
\el
\bpr

By the CSQST property of $\tau_{1}^{{\boldsymbol\alpha}}$, we see that 
\begin{eqnarray*}
	\Phi_{t}^{{\boldsymbol\alpha}}(y) & = & \frac{\mathbb{P}\left(X_{t}^{\alpha}=y\right)-\mathbb{P}\left({\boldsymbol X}_{t}^{{\boldsymbol\alpha}}=(y,1)\right)}{\mathbb{P}(\tau_{1}^{{\boldsymbol\alpha}}>t)}\\
	& = & K'\mu_{t}^{\alpha}(y)-H'\mu^{*}(y)
\end{eqnarray*}
with $K'=1/\mathbb{P}(\tau_{1}^{{\boldsymbol\alpha}}>t)$ and $H'=\mathbb{P}(\tau_{1}^{{\boldsymbol\alpha}}\le t<\tau_{G}^{\alpha})/\mathbb{P}(\tau_{1}^{{\boldsymbol\alpha}}>t)$.
By plugging this equation into the definition of $\tilde{\Phi}_{t}^{\tilde\alpha}$
, we get

\begin{eqnarray}
\frac{\Phi_{t}^{{\boldsymbol\alpha}}(y)\gamma(y)}{\Phi_{t}^{{\boldsymbol\alpha}}\cdot\gamma} & = & \frac{\sum_{x\in A}\alpha(x)\left(K'P_{x,y}^{t}-H'\mu^{*}(y)\right)\gamma(y)}{\Phi_{t}^{{\boldsymbol\alpha}}\cdot\gamma},\label{eq:Fi1}
\end{eqnarray}
by using \eqref{tildePt} and $\nu(x)=\mu^{*}(x)\gamma(x)$,
we get
\begin{eqnarray*}
	\text{(r.h.s. of \ref{eq:Fi1})} & = & K''\sum_{x\in A}\alpha(x)\frac{\gamma(x)}{\gamma(y)}\lambda^{t}\widetilde{P}_{x,y}^{t}\gamma(y)-H\nu(y),
\end{eqnarray*}
where $K''=\frac{K'}{\Phi_{t}^{{\boldsymbol\alpha}}\cdot\gamma}$ and $H=\frac{H'}{\Phi_{t}^{{\boldsymbol\alpha}}\cdot\gamma}$
do not depend on $y$. Since $\tilde{\alpha}(x)=\frac{\alpha(x)\gamma(x)}{\alpha\cdot\gamma}$,
we immediately get the thesis with $K=\lambda^{t}K''\alpha\cdot\gamma$.

\epr

With the help of Lemma \ref{KH}, we can prove the following iteration
formula for the separation $\tilde{s}^{\tilde{\alpha}}(t,y)$:

\bl{sfi}

For any initial measure $\tilde{\alpha}$ on $A$ , any time $t$
and $u$ and $z\in A$, there exist two functions $U$ and $V$, that
depend on $\tilde{\alpha}$ and $t$, but not on $z$ nor on $u$,
such that
\[
\tilde{s}^{\tilde{\alpha}}(t+u,z)=U\tilde{s}^{\tilde{\Phi}_{t}^{\tilde\alpha}}(u,z)+V
\]
\el

\bpr

\begin{eqnarray}
1-\tilde{s}^{\tilde{\Phi}_{t}^{\tilde\alpha}}(u,z) & = & \frac{\sum_{x\in A}\tilde{\Phi}_{t}^{{\alpha}}(x)\widetilde{P}_{x,z}^{u}}{\nu(z)}.\label{eq:comps1}
\end{eqnarray}

By Lemma \ref{KH}, 
\begin{eqnarray*}
	\text{(r.h.s of \ref{eq:comps1})} & = & \frac{\sum_{x\in A}\left(K\tilde{\mu}_{t}^{\tilde{\alpha}}(x)-H\nu(x)\right)\widetilde{P}_{x,z}^{u}}{\nu(z)}\\
	= K\frac{\tilde{\mu}_{t+u}^{\tilde{\alpha}}(z)}{\nu(z)}-H & = & K\left(1-\tilde{s}^{\tilde{\alpha}}(t+u,z)\right)-H
\end{eqnarray*}
and the thesis follows immediately with $U=1/K$ and $V=1-\frac{1+H}{K}$.

\epr

A corollary of this result is the following

\bl{Jf}

For any initial measure $\alpha$ on $A$ , any time $t$ and $u$
and $z\in\mathcal{X}$

\[
J^{\alpha}({t+u},z)=J^{\Phi_{t}^{{\boldsymbol\alpha}}}(u,z)
\]
\el

\bpr

By direct computation, if $z\in A$, 

\begin{eqnarray*}
	J^{\alpha}({t+u},z)=\frac{\tilde{s}^{\tilde{\alpha}}(t+u-1)-\tilde{s}^{\tilde{\alpha}}(t+u)}{\tilde{s}^{\tilde{\alpha}}(t+u-1)-\tilde{s}^{\tilde{\alpha}}(t+u,z)} & = & \frac{\tilde{s}^{\tilde{\Phi}_{t}^{\tilde\alpha}}(u-1)-\tilde{s}^{\tilde{\Phi}_{t}^{\tilde\alpha}}(u)}{\tilde{s}^{\tilde{\Phi}_{t}^{\tilde\alpha}}(u-1)-\tilde{s}^{\tilde{\Phi}_{t}^{\tilde\alpha}}(u,z)}=J^{\Phi_{t}^{{\boldsymbol\alpha}}}(u,z);
\end{eqnarray*}
while, for $z\in G$, we have $J^{\alpha}({t+u},z)=J^{\Phi_{t}^{{\boldsymbol\alpha}}}(u,z)=1$.

\epr

{\bf Proof of Proposition \ref{p_tau1}.}

By the Markov property and Lemma \ref{Jf}, 
$$
\mathbb{P}(\tau_{1}^{{\boldsymbol\alpha}}>t+u)=
$$
\begin{eqnarray*}
	 & = & \sum_{x_{t}\in A}\mathbb{P}({\boldsymbol X}_{t}^{{\boldsymbol\alpha}}=(x_{t},0))
	\sum_{x_{t+1},...,x_{t+u}}\prod_{v=1}^{u}
P(x_{t+v-1},x_{t+v})\Big(1-J^{\a}(t+v,x_{t+v})\Big)\\
	& = & \mathbb{P}(\tau_{1}^{{\boldsymbol\alpha}}>t)\sum_{x_{t}\in A}\Phi_{t}^{{\boldsymbol\alpha}}(x_{t})\sum_{x_{t+1},...,x_{t+u}}\prod_{v=1}^{u}P(x_{t+v-1},x_{t+v})\Big(1-J^{\Phi^{{\boldsymbol\alpha}}_t}(v,x_{t+v})\Big)\\
	& = & \mathbb{P}(\tau_{1}^{{\boldsymbol\alpha}}>t)\mathbb{P}(\tau_{1}^{\Phi_{t}^{{\boldsymbol\alpha}}}>u)
\end{eqnarray*}
\cvd

{\bf Proof of Proposition \ref{p_eph}.}

With the same expansion we write
$$\Phi_{t+u}^{{\boldsymbol\alpha}}(x_{t+u})=$$
\begin{eqnarray*}
& = &\sum_{x_{t}\in A}\frac{\mathbb{P}({\boldsymbol X}_{t}^{{\boldsymbol\alpha}}=(x_{t},0))}{ \mathbb{P}(\tau_{1}^{{\boldsymbol\alpha}}>t+u)}
	\sum_{x_{t+1},...,x_{t+u-1}}\prod_{v=1}^{u}
P(x_{t+v-1},x_{t+v})\Big(1-J^\a(t+v,x_{t+v})\Big)\\
& = &\frac{ \mathbb{P}(\tau_{1}^{{\boldsymbol\alpha}}>t)}{ \mathbb{P}(\tau_{1}^{{\boldsymbol\alpha}}>t+u)}\sum_{x_{t}\in A}
\Phi_{t}^{{\boldsymbol\alpha}}(x_t)
	\sum_{x_{t+1},...,x_{t+u-1}}\prod_{v=1}^{u}
P(x_{t+v-1},x_{t+v})\Big(1-J^{\Phi^{{\boldsymbol\alpha}}_t}(v,x_{t+v})\Big)\\
& = &
\frac{ \mathbb{P}(\tau_{1}^{{\boldsymbol\alpha}}>t)\mathbb{P}(\tau_{1}^{\Phi^{{\boldsymbol\alpha}}_t}>u)}{ \mathbb{P}(\tau_{1}^{{\boldsymbol\alpha}}>t+u)}
\Phi^{\Phi^{{\boldsymbol \a}}_t}_u(x_{t+u})
\end{eqnarray*}
and by Proposition  \ref{p_tau1} we conclude the proof.\cvd
\color{black}
%%%%%%%%%%%%%%%%%%
%%%%
\subsection{Submultiplicativity of $\sup_{\alpha}\mathbb{P}(\tau^\alpha_{*,G}>t)$}
\label{sub}

In this section we prove  Theorem \ref{submult}.

Let $\tau_{*}^{\alpha}$ be a minimal CSQST, and $\tau_{*,G}^{\alpha}=\tau_{*}^{\alpha}\wedge\tau_{G}^{\alpha}$
the associated local relaxation time. We prove that
the function $f(t):=\sup_{\alpha}\mathbb{P}(\tau_{*,G}^{\alpha}>t)$
is submultiplicative, i.e. that for any $t,u>0$, $f(t+u)\le f(t)f(u)$.
This fundamental property implies that $\mathbb{P}(\tau_{*,G}^{\alpha}>t)$
has an exponential bound, allowing in the next section to estimate
the error terms in \eqref{main2_1}.

We start by observing that it is sufficient to study a particular
realization of a minimal CSQST, since 
$$\mathbb{P}(\tau_{*,G}^{\alpha}>t)=\sum_{y\in A}\mathbb{P}\left(X_{t}^{\alpha}=y;\:\tau_{*}^{\alpha}>t\right)=
\sum_{y\in A}  \Big(\mu_{t}^{\alpha}(y)-\mu^{*}(y)\mathbb{P}(\tau_{*}^{{\boldsymbol\alpha}}\le t<\tau_{G}^{\alpha}) \Big)=$$
$$
1-\mu_{t}^{\alpha}(G)-\lambda^{t+\delta_{\alpha}}(1-\tilde{s}^{\tilde{\alpha}}(t))$$
does not depend on the choice of the minimal CSQST. 

Consider a particular realization of the
minimal CSQST, namely, the time $\tau_{1}^{{\boldsymbol\alpha}}$ defined in Section \ref{trp}.
Applying now Proposition \ref{p_tau1} we immediately complete the proof.

%%%%%%%%%%
\subsection{Representation formula for $\t^{\a}_G$ with $\t^\a_*$} 
\label{s_main2}

In this section we prove Theorem \ref{main2}.
Equation (\ref{main2_1}) is an immediate consequence of  the definition of minimal CSQST and of Theorem \ref{main1}.
%%%%%%%%%%%

To prove the final statement of the theorem on the hitting distribution note that for any $y\in G$ we have
 $$
 \P\Big(X^{\a}_{\t^{\a}_G}=y\Big)=\P\Big( \t^{\a}_G< \t^{\a}_*, \; X^{\a}_{\t^{\a}_G}=y\Big)+
\P\Big( \t^{\a}_G> \t^{\a}_*, \; X^{\a}_{\t^{\a}_G}=y\Big) 
$$
The second term in the r.h.s. can be written as
$$
\sum_{t=0}^\infty\sum_{z\in A}\P\Big( \t^{\a}_G> t=\t^{\a}_*, \: X^{\a}_t=z\Big)\P\Big(X^z_{\t^{z}_G}=y\Big)
=\sum_{t=0}^\infty\sum_{z\in A}\m^*(z)\P\Big( \t^{\a}_G> t=\t^{\a}_*\Big)\P\Big(X^z_{\t^{z}_G}=y\Big)$$

$$=\P\Big( \t^{\a}_G>\t^{\a}_*\Big)\sum_{z\in A}\m^*(z)\sum_{u=0}^{\infty}\sum_{w\in A}\P(X^z_u=w, \; \t^z_G=u+1)=
\omega(y)\P\Big( \t^{\a}_G> \t^{\a}_*\Big)
$$
so that (\ref{uscita}) holds.
\bigskip

%%%%%%%%%%%%%%%%%%%%%%%%%%%%%%
\subsection{Under the mean metastability hypothesis}
In this section we prove Theorem  \ref{mainexp}.

  By Theorem \ref{main2} we have
  $$
       \frac{\P\Big( \t^{\a}_{G}> n T \Big)}{\lambda^{nT+\delta_\alpha} } -1= -\tilde s^{\tilde\a}(nT)+
       \frac{\mathbb{P}(\tau_{*,G}^{\alpha}>nT)}{\lambda^{nT+\delta_\alpha}} .
   $$

  %%%%%%%%%%%%%%%%%
  %%
  %%%%%%%%%%%%%%%%%%
 By applying Theorem \ref{submult} and the Markov inequality we have

\be{mexp1}
  \mathbb{P}(\tau_{*,G}^{\alpha}>nT) \le\left(\sup_{\alpha}\mathbb{P}(\tau_{*,G}^{\alpha}>T)\right)^{n}\le \Big(\frac{R}{T}\Big)^{n}.
\ee
We can prove the upper bound:
$$
 \frac{\P\Big( \t^{\a}_{G}> n T \Big)}{\lambda^{nT+\delta_\alpha} } -1\le
 \Big(\frac{R}{T\l^T}\Big)^n \l^{-\d_a}=\frac{e^{-an}}{\a\cdot\g}.
 $$

  As for the lower bound, notice that, by Theorem \ref{main1} and by the minimality of $\tau_{*}^{\alpha}$, $$\mathbb{P}(\tau_{*,G}^{\alpha}=t)\ge\mathbb{P}(\tau_{*}^{\alpha}=t<\tau_{G}^{\alpha})=\lambda^{t+\delta_{\alpha}}\left(\tilde{s}^{\tilde{\alpha}}(t-1)-\tilde{s}^{\tilde{\alpha}}(t)\right),$$
so that
\begin{eqnarray}
	\tilde{s}^{\tilde\alpha}(t) & = & \sum_{u>t}\tilde{s}^{\tilde{\alpha}}(u-1)-\tilde{s}^{\tilde{\alpha}}(u)\le\sum_{u>t}\lambda^{-u-\delta_{\alpha}}\mathbb{P}(\tau_{*,G}^{\alpha}=u)\nonumber\\
	& = & \sum_{u>t}\lambda^{-u-\delta_{\alpha}}\left(\mathbb{P}(\tau_{*,G}^{\alpha}>u-1)-\mathbb{P}(\tau_{*,G}^{\alpha}>u)\right)\nonumber\\
	& = & \lambda^{-t-\delta_{\alpha}}\mathbb{P}(\tau_{*,G}^{\alpha}>t)+\frac{1-\lambda}{\lambda}\sum_{u>t}\lambda^{-u-\delta_{\alpha}}\mathbb{P}(\tau_{*,G}^{\alpha}>u).\label{mexp2}
\end{eqnarray}
Thus, the total error in \eqref{main2_1} can be
bounded as

\begin{eqnarray}
\mathbb{P}(\tau_{*,G}^{\alpha}>t)-\lambda^{t+\delta_{\alpha}}\tilde{s}^{\tilde{\alpha}}(t) & \ge & -\frac{1-\lambda}{\lambda} \lambda^t \sum_{u>t}\lambda^{-u}\mathbb{P}(\tau_{*,G}^{\alpha}>u).\label{eq:Lerr}
\end{eqnarray}

Again by Markov inequality and submultiplicativity,
\begin{eqnarray*}
  \sum_{u>nT}\lambda^{-u}\mathbb{P}(\tau_{*,G}^{\alpha}>u)
  \le
  T  \sum_{k\ge n}\lambda^{-(k+1) T}\Big(\frac{R}{T}\Big)^k
  =  \frac{T \lambda^{-T+1}}{1-\frac{R}{T\lambda^T }} \left(\frac{R}{T\lambda^T }\right)^n, 
\end{eqnarray*}
where we used the fact that the sum is convergent since $\frac{R}{T\l^T}<1$ by hypothesis. 

Thus, we obtain the lower bound
$$
 \frac{\P\Big( \t^{\a}_{G}> n T \Big)}{\lambda^{nT+\delta_\alpha} } -1\ge
 -\frac{1-\lambda}{\lambda} \lambda^{-\d_a} \sum_{u>nT}\lambda^{-u}\mathbb{P}(\tau_{*,G}^{\alpha}>u)\ge
$$
$$
-\l^{-\d_a} \frac{\l^{-T}}{1-e^{-a}}e^{-an}
$$
and the thesis immediately follows.

%%%%%%%%%%%%%%%%%%%%%%%%%%%%%%%%%%%

%%%%%%%%%%%%
\section{An example: the rim}
\label{ex}

As explained  in the introduction, metastability is  associated
to the existence of two asymptotically-separated time-scales: a ``short''
time-scale in which the system relaxes to a sort of ``apparent equilibrium''
and a ``long'' time-scale that characterizes the arrival to the
invariant measure. We introduce here a simple model, inspired by a
similar example introduced in \cite{DF90} to show how the
Conditionally Strong Time language can be used to formalize this picture.

In this model, each term in the representation formula can be computed,
so that we can illustrate the meaning of the terms of the representation formula \eqref{main2_1}
in an explicit case.

\medskip{}

Let $n\in\mathbb{N}$ be an integer parameter. The state-space of
the model is $\mathbb{T}_{n}\cup G$, where $\mathbb{T}_{n}$ denotes the 1-dimensional discrete 
torus of lenght $4^{n}$(labeled from $0$ to $4^{n}-1$) and $G$
is a single absorbing state. The graph is illustrated in fig. 1.

\begin{figure}[h]
\centering
\includegraphics[width=0.5\linewidth]{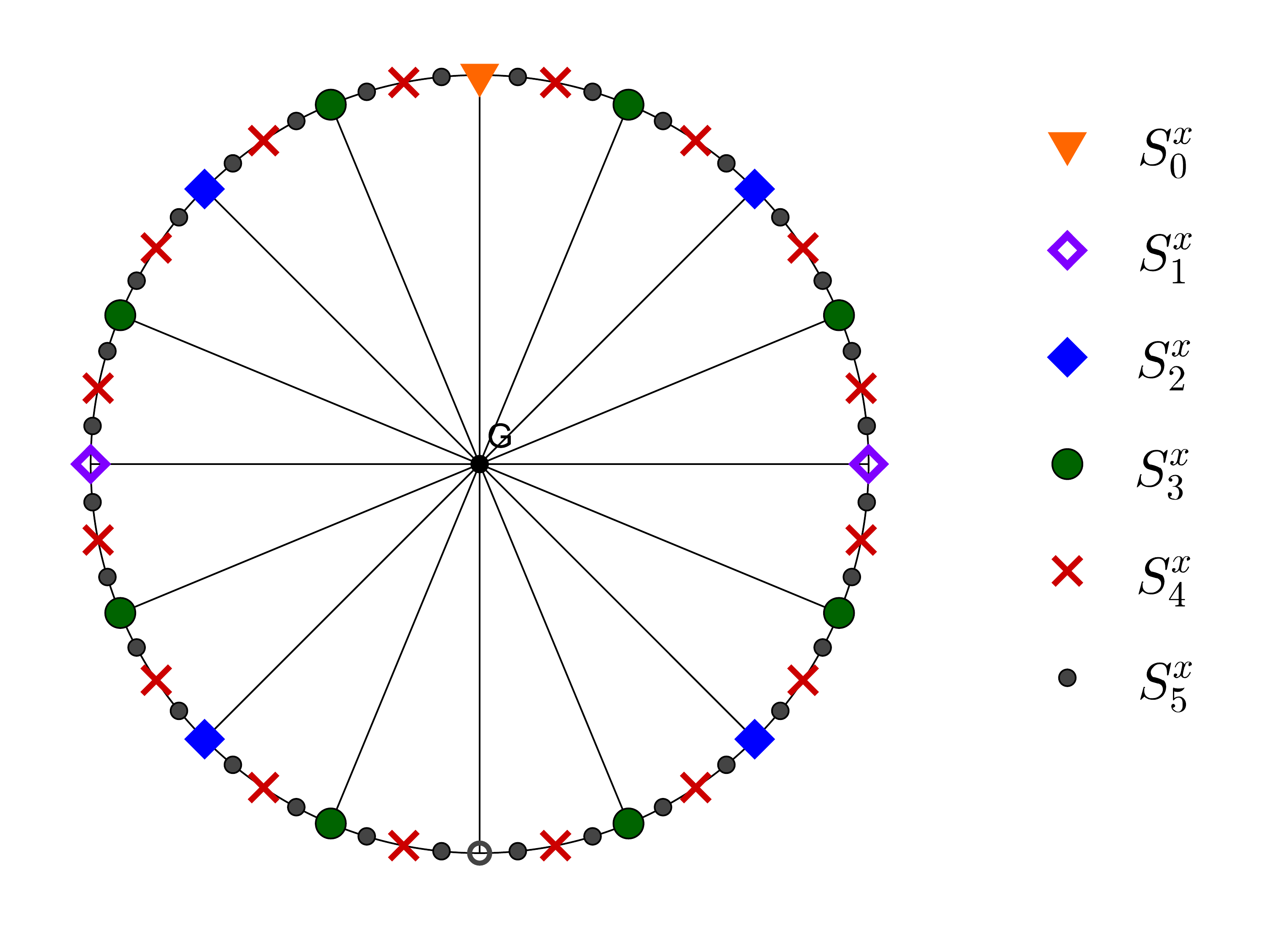}
\caption{The graph and the subsets $S_{k}^{x}$ when $n=3$}
\label{fig:grafo}
\end{figure}

All transition probabilities will be invariant under rotations by
multiples of $4$ , so that,
 according to the heuristic definition
of metastability given above, we should compare the time to diffuse
onto the ring to the time needed to take one of the spokes. If the
former  time is much shorter than the latter, the system somehow
thermalizes before undergoing the transition to equilibrium; if not,
we cannot talk about metastability even when the arrival time to $G$
is exponentially-distributed.

The transition probabilities, as we will show later on, are chosen
to keep as simple as possible the construction of the CSQST. On the
same graph, all choices with similar symmetries would allow the construction
of the CSQST, but the construction would not be as simple. 

Let $\mathbb{T}_{n}^{0}$ be the subset of the multiples of $4$ in $\mathbb{T}_{n}$
, $\mathbb{T}_{n}^{1}$ denote the subset of the odd numbers in $\mathbb{T}_{n}$ and
$\mathbb{T}_{n}^{2}=\mathbb{T}_{n}\setminus \Big(\mathbb{T}_{n}^{0}\cup \mathbb{T}_{n}^{1}\Big)$ denote the
remaining subset.

Let $\lambda\in (0,1)$ be a real parameter that will correspond to the largest eigenvalue of the matrix $[P]_{A}$.
The non-null elements of the transition matrix $P$ are:

If $x\in \mathbb{T}_{n}^{0}$ (that is, a multiple of $4$),

\[
P_{x,y}=\begin{cases}
\frac{\lambda}{2} & \text{\text{ if }}y=x\\
\frac{\lambda^{2}(2-\lambda)}{32-32\lambda+4\lambda^{2}} & \text{\text{ if }}y\in \mathbb{T}_{n}\text{ and }|x-y|=1\\
\frac{8-12\lambda+4\lambda^{2}}{8-8\lambda+\lambda^{2}} & \text{\text{ if }}y=G
\end{cases}
\]

If $x\in \mathbb{T}_{n}^{1}$ (that is, an odd number),

\[
P_{x,y}=\begin{cases}
\frac{\lambda}{2} & \text{\text{ if }}y=x\\
\frac{8-8\lambda+\lambda^{2}}{4(2-\lambda)} & \text{\text{ if }}y\in \mathbb{T}_{n}^{0}\text{ and }|x-y|=1\\
\frac{\lambda^{2}}{4(2-\lambda)} & \text{\text{ if }}y\in \mathbb{T}_{n}^{2}\text{ and }|x-y|=1
\end{cases}
\]

If $x\in \mathbb{T}_{n}^{2}$ (that is, an even number but not a multiple of
$4$),

\[
P_{x,y}=\begin{cases}
\frac{\lambda}{2} & \text{\text{ if }}y=x\\
\frac{2-\lambda}{4} & \text{\text{ if }}y\in \mathbb{T}_{n}\text{ and }|x-y|=1
\end{cases}
\]

Moreover, $P_{G,G}=1$.

\medskip{}

The Perron-Frobenius theorem and a direct computation allows to prove the following:

\bp{ex1}
$\lambda$ is the largest eigenvalue of the sub-markovian
matrix $[P]_{A}$, associated to the left eigenvector (normalized
to $1$)

\[
\mu_{x}^{*}=\begin{cases}
4^{-n}\frac{8-8\lambda+\lambda^{2}}{2-\lambda} & \text{ if }x\in \mathbb{T}_{n}^{0}\\
4^{-n}\lambda & \text{ if }x\in \mathbb{T}_{n}^{1}\\
4^{-n}\frac{\lambda^{2}}{(2-\lambda)} & \text{ if }x\in \mathbb{T}_{n}^{2}
\end{cases}
\]
 and to the right eigenvector

\[
\gamma_{x}=\begin{cases}
\frac{2-\lambda}{8-8\lambda+\lambda^{2}} & \text{ if }x\in \mathbb{T}_{n}^{0}\\
\frac{1}{\lambda} & \text{ if }x\in \mathbb{T}_{n}^{1}\\
\frac{2-\lambda}{\lambda^{2}} & \text{ if }x\in \mathbb{T}_{n}^{2}
\end{cases}
\]
with normalization such that $\sum_{x}\gamma_{x}\mu_{x}^{*}=1$.
\ep
\bigskip

\br{R1}
If $x\in \mathbb{T}_{n}^{1}$ and $y=x\pm1$, then $P_{x,y}=4^{n-1}\mu_{y}^{*}$,
while $P_{x,x}=2\cdot4^{n-1}\mu_{x}^{*}$. Therefore, starting from
the uniform distribution $\alpha$ on $\mathbb{T}_{n}^{1}$:
$$
\mu^\alpha_1(y)\equiv P_{\alpha, y}=\sum_{x\in \mathbb{T}_{n}^{1}}2\cdot {4^{-n}}P(x,y)=\mu^*_y
$$
\er
\medskip{}

In order to construct the CSQST, we define a family of sets $S_{k}^{x}$
recursively:

Let $x$ be a starting configuration, for $k\in\left\{ 1,\ldots,2n-1\right\} $

\begin{eqnarray*}
\tau^{x}(0) & := & \tau_{\mathbb{T}_{n}^{0}}^{x}\\
S_{0}^{x} & := & \left\{ X_{\tau^{x}(0)}^{x}\right\} \\
S_{k}^{x} & := & \left\{ y\in \mathbb{T}_{n}\:;\:y\pm2^{2n-k-1}\in S_{k-1}^{x}\right\} \\
\tau^{x}(k) & := & \inf\left\{ t>\tau^{x}(k-1)\:;\:X_{t}^{x}\in S_{k}^{x}\right\} 
\end{eqnarray*}

where the symbol $\pm$ denotes the sum/difference modulo $4^{n}$ (see Figure 1).
\bigskip

\br{R2} 
\ 
\begin{itemize}
\item[1)] For every $k\in \{2,\dots,2n-2\}$, between each two consecutive elements of $S_{k-1}^{x}$we
put two elements of $S_{k}^{x}$.
\item[2)] $|S_{k}^{x}|=2^{k}$.
\item[3)] The sets $S_{k}^{x}$ are stochastic only because $S_{0}^{x}$ is stochastic.
As we will see in what follows, due to the symmetry of the model, we 
are mainly interested in the case
 $x=0$. In this case $S_{0}^{x}=0$.
\item[4)] 
 To each element $y$ of $S_{k}^{x}$ is associated
a unique ``parent'' $g(y)$ in $S_{k-1}^{x}$ such that $|y-g(y)|=2^{2n-k-1}$. 
Each parent has two offsprings and
 for every $k\ge1,$ if $X_{\tau^{x}(k)}^{x}=y$, then $X_{\tau^{x}(k-1)}^{x}=g(y)$.
 
 \end{itemize}
 \er
 With these definitions we can state our main result on this model.
\bigskip

\bt{thex}
The time
\[
\tau_{\star}^{x}:=\tau^{x}(2n-1)+1
\]
is a CSQST.
\et
\bigskip

{\bf Proof}:

We first consider the case $x=0$ and
we start by proving inductively that for each $k\in\left\{ 0,\ldots,2n-1\right\} $,
$X_{\tau^{0}(k)}^{x}$is independent of $\tau^{0}(k)$ and uniformly
distributed on $S_{k}^{0}$: for $y\in S_{k}^{0}$, 
\begin{equation}
\mathbb{P}\left(X_{t}^{0}=y,\:\tau^{0}(k)=t\right)=2^{-k}\mathbb{P}\left(\tau^{0}(k)=t\right).\label{eq:unifA}
\end{equation}
Indeed, for $k=0$, $S_{0}^{0}=0$ is a singleton and there is nothing
to prove. For $k\ge1$, 
\[
\mathbb{P}\left(X_{t}^{0}=y,\:\tau^{0}(k)=t\right)=\sum_{s<t}\mathbb{P}\left(X_{t}^{0}=y,\:\tau^{0}(k)=t,\:X_{s}^{0}=g(y),\:\tau^{0}(k-1)=s\right)
\]
by symmetry we have
\[
\mathbb{P}\left(X_{t}^{0}=y,\:\tau^{0}(k)=t\right)=\frac{1}{2}\sum_{s<t}\mathbb{P}\left(\tau^{0}(k)=t,\:X_{s}^{0}=g(y),\:\tau^{0}(k-1)=s\right)
\]
and by the inductive hypothesis
we get 
\[
=\frac{1}{2}\sum_{s<t}\mathbb{P}\left(\tau^{0}(k)=t\big|X_{s}^{0}=g(y),\:\tau^{0}(k-1)=s\right) 2^{-k+1}
\mathbb{P}\left(\tau^{0}(k-1)=s\right)
\]
By symmetry we can ignore the conditioning $X_{s}^{0}=g(y)$ obtaining 
\[
=\frac{1}{2}\sum_{s<t}\mathbb{P}\left(\tau^{0}(k)=t,\:\tau^{0}(k-1)=s\right) 2^{-k+1}=2^{-k}
\mathbb{P}\left(\tau^{0}(k)=t\right)
\]
Now, we observe that $S_{2n-1}^{0}$is the set $\mathbb{T}_{n}^{1}$ of the
odd numbers, from which it is not allowed to move to $G$, and by recalling Remark \ref{R1}
we immediately obtain for any $y\in \mathbb{T}_{n}$:
$$
\mathbb{P}\left(X_{t}^{0}=y,\:\tau_{\star}^{0}=t\right)=\mu_{y}^{*}\mathbb{P}\left(\tau_{\star}^{0}=t\right)
$$
Since in our example $\tau_{\star}^{0}=t$ entails $\tau_{G}^{0}>t$,
we can condition the last two formulae and see that $\tau_{\star}^{0}$
is a CSQST.
\bigskip

If we now consider an arbitrary starting point $x$, due to the definition of $S^x_0$ and $\tau^x_0$,
 this is equivalent to start from some point in $\mathbb{T}^0_n$, defined below as $S^x_0$, depending on $x$. 
 As noted above  
this point $S^x_0$ is random if $x\not\in \mathbb{T}^0_n$. However, due to the symmetry of the model, every starting point in $\mathbb{T}^0_n$
is equivalent to $0$. This concludes the proof of the theorem.

\medskip{}

We want to discuss now the application of our representation formula to this particular model.

Let us call $\mathrm {Opp}(x):=S_{0}^{x}+2^{2n-1}$ the point opposite to $S_{0}^{x}$.
In order to reach this point, the process has to visit every set $S_{k}^{x}$,
so that we get the estimates 
$$\mathbb{P}\left(\tau_{\star}^{x}>t\right)\le\mathbb{P}\left(\tau_{\mathrm {Opp}(x)}^{x}\ge t\right)\quad
\hbox{
and }\quad\mathbb{P}\left(\tau_{\star,G}^{x}>t\right)\le\mathbb{P}\left(\tau_{\mathrm {Opp}(x),G}^{x}\ge t\right).
$$
Standard diffusive bounds show that $\mathbb{P}\left(\tau_{\mathrm {Opp}(x)}^{x}\le 4^{2n-1}\right)$ is
larger than a constant $c\in(0,1]$. By dividing the time $t$ into
intervals of length $4^{2n-1}$, we get 
\begin{equation}
\mathbb{P}\left(\tau_{\mathrm {Opp}(x),G}^{x}\ge t\right)\le\left(\sup_{y\in \mathbb{T}^{n}}\mathbb{P}\left(\tau_{\mathrm {Opp}(x),G}^{x}\ge 4^{2n-1}\right)\right)^{4^{-2n+1}t}\le\left(1-c\right)^{4^{-2n+1}t}\label{eq:es3}
\end{equation}
In order to compute the other terms in the representation formula,
we consider the local chain on $\mathbb{T}_{n}$ 

\[
\widetilde{P}_{x,y}=\frac{\gamma(y)}{\gamma(x)}\frac{P_{x,y}}{\lambda}=\begin{cases}
\frac{1}{2} & \text{ if }x=y\\
\frac{1}{4} & \text{ if }|x-y|=1
\end{cases}
\]

Clearly, the local process $\widetilde{P}_{x,y}^{t}$ is a lazy random
walk on the ring. It is easy to control the convergence to equilibrium for this process in separation distance.
Indeed,  by standard diffusive estimates we get that 
$\tilde{s}^{x}(4^{2n-1})\le b$
for some constant $b\in (0,1)$.

Since $\sup_{x\in \mathbb{T}_{n}}\tilde{s}^{x}(t)$ is submultiplicative, we
get

\begin{equation}
\widetilde{s}^{x}(t)\le b^{4^{-2n+1}t}.\label{eq:es2}
\end{equation}

From these estimates, we see that the error terms in \eqref{main2_1} can be estimated by
$\gamma_x \lambda^t b^{4^{-2n+1}t} + \left(1-c\right)^{4^{-2n+1}t}$
and, when the time needed to diffuse onto the ring is smaller than the mean time $1/(1-\lambda)$ to reach $G$, they decay faster than the leading term $\gamma_x \lambda^t$.

It is useful to compare the estimate given by theorem \ref{main2}
with a direct computation of $\mathbb{P}\left(\tau_{G}^{x}>t\right)$.
To this end, let us introduce the projection operator $p:\:\mathbb{T}_{n}\cup G\rightarrow\left\{ 0,1,2,G\right\} $
defined by $p(G)=G$ and $x\in \mathbb{T}_{n}^{p(x)}$. We notice that the
projection $p(P^{t})$ is itself a Markov process with transition
matrix

$$\bar{P}=\left(\begin{array}{cccc}
P_{0,0} & 2P_{0,1} & 0 & P_{0,G}\\
P_{1,0} & P_{1,1} & P_{1,2} & 0\\
0 & 2P_{2,1} & P_{2,2} & 0\\
0 & 0 & 0 & 1
\end{array}\right)$$

this means that $P_{i,j}=P_{x,y}$ with $x\in \mathbb{T}^i_n,\, y\in \mathbb{T}^j_n$ for any $i,j\in \{0,\, 1,\, 2,\, G\}$.
Thus, $\mathbb{P}\left(\tau_{G}^{x}>t\right)=\bar{\mathbb{P}}\left(\tau_{G}^{p(x)}>t\right)$, where $\bar{\mathbb{P}}$
denotes the probability for the Markov chain with transition matrix $\bar P$. 

The largest eigenvalue of the restricted matrix is again $\lambda$ and
the quasi-stationary measure is the projection of $\mu^{*}$:

\[
\begin{cases}
\bar{\mu}_{0}^{*}= & \frac{1}{4}\frac{8-8\lambda+\lambda^{2}}{2-\lambda}\\
\bar{\mu}_{1}^{*}= & \frac{\lambda}{2}\\
\bar{\mu}_{2}^{*}= & \frac{\lambda^{2}}{4(2-\lambda)}
\end{cases}
\]

\br{R3}
 $\bar{\mu}_{i}^{*}=\bar{P}_{1,i}$. Hence, when starting
from $1$, the projected chain reaches equilibrium at time $1$.
\er
Thus, 
\[
\bar{\mathbb{P}}\left(\tau_{G}^{1}>t\right)=\lambda^{t-1}= \gamma_1 \lambda^t
\]
 
\begin{eqnarray*}
\bar{\mathbb{P}}\left(\tau_{G}^{2}>t\right)=
\bar P_{2,2}^t+
\sum_{s=0}^{t-1}\bar{P}_{2,2}^{s}\bar{P}_{2,1}\lambda^{t-s-2}=
\left(\frac{\lambda}{2}\right)^t+ \sum_{s=0}^{t-1}\left(\frac{\lambda}{2}\right)^s\frac{2-\lambda}{2}
\lambda ^{t-s-2}=\\
\left(\frac{\lambda}{2}\right)^t+ \frac{2-\lambda}{2}\lambda ^{t-2}2(1-2^{-t})=
\gamma_2\lambda^t\big(1-2^{-t}(1-\frac{1}{\gamma_2})\big),
\end{eqnarray*}

\[
\bar{\mathbb{P}}\left(\tau_{G}^{0}>t\right)=
\bar{P}_{0,0}^{t}+\sum_{s=0}^{t-1}\bar{P}_{0,0}^{s}\bar{P}_{0,1}\lambda^{t-s-2}=
\left(\frac{\lambda}{2}\right)^t+\frac{\lambda^2}{4}\gamma_0 4\lambda^{t-2}(1-2^{-t})=
\gamma_0\lambda^t\big(1-2^{-t}(1-\frac{1}{\gamma_0})\big).
\]

We see that, due to the symmetry of this system, the distribution of $\tau^x_G$ can be approximated with an exponential distribution much before the diffusive time on the ring.
In other words, the hitting time has an exponential behavior even before the metastable time.

{\bf Acknowledgments: }

We thank Amine Asselah, Nils Berglund,  Pietro Caputo, Frank den Hollander, Roberto Fernandez and Alexandre Gaudilli\`ere for many fruitful  discussions.
 This work was partially supported by the A*MIDEX project (n. ANR-11-IDEX-0001-02) funded by the  ``Investissements d'Avenir" French Government program, managed by the French National Research Agency (ANR).
E.S. has been supported by the PRIN 20155PAWZB ÒLarge Scale Random StructuresÓ.
%%%%%%%%%%%%%
%%%%%%%%%%%%%

\end{document}